\documentclass[12pt]{article}  

\usepackage{hyperref}

\usepackage{amsfonts}

\usepackage{enumerate}

\let\oldmarginpar\marginpar
\renewcommand\marginpar[1]{${}^\clubsuit$\oldmarginpar[\raggedleft\scriptsize\sf #1]{\raggedright\scriptsize\sf #1}}
\newtheorem{theorem}{Theorem}

\newtheorem{proposition}[theorem]{Proposition}

\newtheorem{corollary}[theorem]{Corollary}
\newtheorem{remark}[theorem]{Remark}

\newcommand{\proof}{ \noindent{\it Proof:\ \ }}
\renewcommand{\P}{\Bbb{P\,}}
\def\qed{\ifhmode\unskip\nobreak\fi\ifmmode\ifinner\else\hskip5 pt
\fi\fi\hbox{\hskip5 pt \vrule width4 pt height6 pt depth1.5 pt
\hskip 1pt }}
\newcommand{\po}{{\hspace*{-1ex}}{\bf .  }}

\newcommand{\cref}[1]{Corollary~\ref{#1}}

\newcommand{\eqref}[1]{(\ref{#1})}

\newcommand{\pref}[1]{Proposition~\ref{#1}}
\newcommand{\rref}[1]{Remark~\ref{#1}}

\newcommand{\tref}[1]{Theorem~\ref{#1}}

\def\Bbb#1{\mathbb#1}
\newcommand{\R}{\Bbb{R}}

\newcommand{\Sp}{\Bbb{S}}

\newcommand{\C}{\Bbb{C\,}}
\newcommand{\CP}{\Bbb{C}\Bbb{P}\,}

\newcommand{\imc}[4]{#1\colon #2^{#3} \rightarrow \C^{#3+#4}}

\newcommand{\nn}[1]{\|#1\|^{-2}}
\newcommand{\n}[1]{\|#1\|^{2}}
\newcommand{\sff}{second fundamental form }

\newcommand{\nape}{\nabla^\perp}

\renewcommand{\a}{\alpha}

\newcommand{\ra}{\rangle}
\newcommand{\la}{\langle}

\newcommand{\spa}{\mbox{span}}

\newcommand{\cl}{{\cal L}}

\newcommand{\hisim}{holomorphic isometric immersion }
\newcommand{\h}{holomorphic}
\newcommand{\na}{\nabla}
\newcommand{\fp}{{f^\perp}}
\newcommand{\ft}{{f^\top}}
\newcommand{\hf}{{\hat f}}
\newcommand{\hfp}{{\hf^\perp}}
\newcommand{\pft}{{\a_\ft}}
\newcommand{\gp}{{g^\perp}}

\newcommand{\hm}{{\hat M^{n-\nu}}}
\newcommand{\K}{K\"ahler }
\newcommand{\rank}{\mbox{rank }}

\begin{document}

\title{The holomorphic Gauss Parametrization  
}
\author{Marcos Dajczer \& Luis A. Florit}

\date{}
\maketitle

\begin{abstract}  
We give a local parametric description of all complex hypersurfaces
in $\C^{n+1}$ and in complex projective space $\CP^{n+1}$ with
constant index of relative nullity, together with applications.
This is a complex analogue to the parametrization for real
hypersurfaces in Euclidean space known as the Gauss parametrization.

\end{abstract}

\section{Introduction}  

Let $M^n$ be a complete connected complex immersed hypersurface of
$\C^{n+1}$ whose index of relative nullity, that is, the
dimension of the kernel of its second fundamental form, satisfies
$\nu\geq n-1$ everywhere. Equivalently, its Gauss map
$\varphi\colon\,M^n\to\CP^n$ that assigns to each point in $M^n$ its
normal complex line in $\C^{n+1}$, satisfies $\rank d\varphi\leq 1$.
Then, it was shown by Abe (\cite{abe}) that the hypersurface must be an
$(n-1)$-cylinder.

The situation is even more restrictive for a complete hypersurface $M^n$
of the complex projective space $\CP^{n+1}$. From a general result also
due to Abe (\cite{abe2}) it follows that if $\nu>0$ then $M^n$ must be a
totally geodesically embedded $\CP^n \subset \CP^{n+1}$.

Naturally, the situation is quite different in the local case. In fact,
for any integer $\nu_0>0$ there are plenty of local hypersurfaces $M^n$
in $\C^{n+1}$ and $\CP^{n+1}$ with constant index of relative nullity
$\nu=\nu_0$ that are neither part of cylinders in $\C^{n+1}$ or totally
geodesic in $\CP^{n+1}$.

Our main goal in this note is to give a parametric description of all
complex hypersurfaces in $\C^{n+1}$ and $\CP^{n+1}$ with constant
$\nu>0$. As a consequence, the above global results will be immediate
corollaries of our local construction achieved by imposing on the
hypersurfaces the absence of singularities.

So far everything just said is an analogue to what happens for real
hypersurfaces in Euclidean space $\R^{n+1}$ and the round sphere
$\Sp^{n+1}$. The now called {\it Gauss parametrization} was introduced by
Sbrana (\cite{sb}) as a tool to classify the locally isometrically
deformable Euclidean hypersurfaces. In recent years, it has proved to be
quite a useful tool, giving rise to several applications; see
\cite{df}, \cite{dft}, \cite{dg}, \cite{dg2}, \cite{dg3}, \cite{dk} and
\cite{dt2}.

The parametrization of hypersurfaces in $\C^{n+1}$ we give here works
similarly for hypersurfaces in $\R^{n+1}$ and provides an equivalent form
of the Gauss parametrization given in~\cite{dg}. Our parametrization for
$\CP^{n+1}$ is a perfect analogue of the Gauss parametrization in the
sphere $\Sp^{n+1}$. We point out that the holomorphicity hypothesis is
redundant for submanifolds in $\CP^{N}$ with $\nu>0$ (\cite{dr}).

\section{The parametrization in $\C^{\!N}$ } 

Let $\imc f M n p$ be a \hisim of a \K Riemannian manifold
$M^n=(M^n,\la\ ,\ \ra)$ with Levi-Civita connection $\na$, normal
connection $\nape$, and \sff $\a:\,TM\oplus TM\to T^\perp M$. We
denote by $J$ the complex structures of both $M^n$ and $\C^{n+1}$.

Recall that the {\it relative nullity} subspace $\Delta(x)$ of $f$
at $x\in M^n$ is given by
$$
\Delta(x)=\{Y\in T_xM: \a(Y,Z) = 0,\ \ \forall\, Z\in T_xM\}.
$$
Since $f$ is holomorphic, $\Delta(x)\subseteq T_xM$ is a complex subspace
whose (complex) dimension $\nu_f(x)$ is called the {\it index of relative
nullity} of $f$ at $x$. Along the open dense subset $M_0\subseteq M^n$
where $\nu$ is constant, it is well known that $\Delta$ is a smooth
integrable distribution whose leaves are totally geodesic in both $M^n$
and $\C^{n+p}$. Locally on a saturated open subset $U\subseteq M_0$, the
space of leaves of this distribution, that we denote by $\hat
M=U/\Delta$, is naturally a complex manifold of dimension $n\!-\!\nu$
whose projection $\pi\colon U\to \hm$ is holomorphic. This space can be
naturally identified with a complex submanifold of $U$ of dimension
$n-\nu$ transversal to the leaves of relative nullity. We point out that
the space of leaves is well defined globally on $M_0$ if $M^n$ is
complete; see~\cite{dg}.

Let $i\colon\,\C_*^N\to\C_*^N$ be the inversion given by $i(z)=z/\|z\|^2$.
Its differential is
$$
di_z(v)= \frac{1}{\|z\|^2}R_zv,
$$
where
$$
R_zv=v-2\la v,z\ra i(z),
$$
stands for the reflection in the $z$ direction. Notice that $R_z$
satisfies
$$
J R_z = R_{Jz} J
$$
where $J$ is the complex structure in $\C^N$.

Let $\imc f M n p$ be a \hisim of a \K  manifold.
Decompose the position vector of the immersion as
\begin{equation}\label{e:des}
f = \fp + \ft,
\end{equation}
according to the orthogonal holomorphic bundle decomposition
$$
\C^{n+p} \cong T_x\C^{n+p} = T_x^\perp M \oplus T_xM,
$$
for each $x\in M^n$.

If there is an open subset $U\subseteq M^n$ for which the
position vector $f$ is tangent to $f$, then by analyticity $f$ is
everywhere tangent and, regarded as a tangent vector field, must belong
to the relative nullity. Hence, $f$ must be a (complex) cone through the
origin. By means of a generic translation $f+p_0$, we assume from now
on that the position vector is not tangent on an open dense subset of
$M^n$, that we continue calling $M^n$.

Differentiating (\ref{e:des}) and taking
normal components, we get
$\pft:=\a(\cdot , \ft) = -\nape \fp$. Hence,
\begin{equation}\label{e:dfp}
d \fp = -A_\fp - \pft,
\end{equation}
where $A_\delta=A^f_\delta$ denotes the real shape operator of $f$ in the direction
$\delta$. Since $f$ is \h, we have that $A_\fp J=-JA_\fp$ and
$\pft J=J\pft$, and thus
\begin{equation}\label{e:hol}
d\fp J = J \left( A_\fp - \pft\right) = -J(d\fp + 2\pft).
\end{equation}
Setting
\begin{equation}\label{e:defg}
g=i(\fp),
\end{equation}
we have
\begin{equation}\label{e:dg}
dg=di_{\fp}\circ d\fp = \frac{1}{\|\fp\|^2}R_\fp d\fp =\n g R_g d\fp,
\end{equation}
since $R_\fp=R_g$. Therefore,
\begin{eqnarray*}
\nn g dgJ\!\! &=&\!\! R_g d\fp J = -R_{g}J(d\fp+2\pft)\\
\!\!&=&\!\! -JR_{Jg}(d\fp+2\pft)\\
\!\!&=&\!\! -J(R_gd\fp + 2 \Pi_\gp(\pft))\\
\!\!&=&\!\! -\nn g Jdg - 2J\Pi_\gp(\pft)
\end{eqnarray*}
where $\Pi_\gp \colon T^\perp M \to  T^\perp M \cap (\spa_\C\{g\})^\perp$
is the orthogonal projection.

In particular, if the codimension is $p=1$ we conclude that
$$
dgJ=-Jdg,
$$
that is, $g$ is anti-holomorphic. Moreover, observe that since $f$ is \h,
by \eqref{e:dfp} and \eqref{e:dg} we have that $\ker dg =
\Delta$ at each point. In other words, $g$ is constant along the leaves of relative
nullity of $f$ and, locally, there is an anti-holomorphic immersion
$\hf\colon \hm \to \C^{n+1}$ such that $g=\hf \circ \pi$. We will always
consider on $\hm$ the \K metric induced by $\hf$.

\vspace{2ex}

We have proved:
\begin{proposition}\po\label{p:antiholo}
Let $\imc f M n 1$ be a \hisim of a \K Riemannian manifold. Then, on
the open dense subset where the position vector $f$ is not tangent, the
map $g=i(\fp)$ is anti-holomorphic. Moreover, locally along the open
dense subset which also has constant index of relative nullity $\nu$,
there is an
anti-\hisim $\hf\colon \hm \to \C^{n+1}$ such that $\hf \circ \pi=g$.
\end{proposition}

Observe that, as a consequence,
the Gauss map $N\colon M^n\to \CP^n$ of $f$ given by
$$
N(x)=\spa_\C\{\hf^\perp(\pi(x))\}
$$
is anti-holomorphic; see \cite{nom}.
\vspace{2ex}

Our purpose now is to describe $f$ locally by means of the geometry of
$\hf$.

\begin{theorem}\po\label{t:param}
Let $\hf\colon\hm \to \C^{n+1}$ be an anti-\hisim of a \K manifold with
$\nu_\hf=0$ whose position vector is never tangent. Let $\cl$
be the holomorphic vector subbundle given by
\begin{equation}\label{e:h}
\cl = {\rm span}_\C\{\hf^\perp\}^\perp\subset T^\perp_\hf\hat M.
\end{equation}
Then, the map $f\colon  \cl \to \C^{n+1}$ defined as
\begin{equation}\label{e:param}
f(\xi) = i(\hf^\perp(x)) + \xi, \ \ \ \xi \in \cl(x),
\end{equation}
parametrizes, at regular points, a \h \ \K hypersurface with constant
index of relative nullity $\nu_f=\nu$.
Conversely, any such hypersurface can be parametrized this way.
\end{theorem}

\proof
For the direct statement, observe first that $f$ is holomorphic by
\pref{p:antiholo}. Moreover, since $\la f , \hf\circ\pi\ra = 1$, we have
that
$$
0=\la df,\hf\circ\pi\ra + \la f,d\hf\circ\pi\ra = \la df,\hf\circ\pi\ra,
$$
that is, $\hf\circ\pi$ is
normal to $f$. From the definition, it is clear that the fibers of $\cl$
are contained in the relative nullity of $f$, and they must coincide
since $\hf$ is never tangent.

For the converse, we follow the arguments before \pref{p:antiholo}
writing
$$
f=i(g) + \ft\renewcommand{\P}{\Bbb{P\,}}
$$
and $\hf\circ \pi=g$, where $g=i(f^\perp)$. Since $g$ is normal to $f$ and
$\ker dg = \Delta$, by dimension reasons we conclude that the leaf of $\Delta$
through $x$ is simply (contained in) a translation of $\cl(\pi(x))$ defined by
\eqref{e:h}. Therefore, we set
$$
\ft = h\circ \pi + \xi
$$
where $h\in\cl^\perp$ and $\xi(x)\in \cl(\pi(x))$. Again by dimension reasons,
$\xi$ parametrizes each leaf of $\cl$ when $x$ moves along a leaf of relative
nullity. Now, differentiating $\la f,g\ra=1$, we obtain $\la f,dg\ra=0$. It
follows that
\begin{equation}\label{e:hp}
i(\hf)+h\in T^\perp\hat M.
\end{equation}
By \eqref{e:h} and \eqref{e:hp} we
have that $i(\hf)+h$ and $\hfp$ are linearly dependent, say,
$i(\hf)+h=\lambda\hfp$
Since $h$ is tangent to $f$, taking the inner product with $\hf$ yields
$i(\hf)+h=i(\hfp)$, as we wanted to prove.
\qed

\begin{remark}\po\label{r:invol}
{\rm
Let ${\cal H}_+$ and ${\cal H}_{-}$ denote the sets of hypersurfaces in
$\C^N$ without relative nullity and whose position vectors are never
tangent, that are holomorphic and anti-holomorphic, respectively. Since
the roles of holomorphic and anti-holomorphic submanifolds can be
reversed in the above arguments, the map defined on
${\cal H}_+\cup{\cal H}_{-}$ given by
$$
f\mapsto f^*=i(f^\perp)
$$
is a bijection that swaps ${\cal H}_+$ with ${\cal H}_{-}$ such that
$(f^*)^*=f$. In the case of holomorphic curves it was shown in \cite{dt}
that this map is conformal.
}\end{remark}

We now compute the singular set and second fundamental form
of the submanifold using the parametrization \eqref{e:param}, the latter
being completely determined by $A_\hf$ by the holomorphicity of $f$.
Let $P\colon\,T\hat M  \to\Delta^\perp$ where
$$
\Delta^\perp=(\cl\oplus\spa_\C\{\hf\})^\perp
=(T\hat M\oplus {\rm span}_\C\{\hf\})\cap
({\rm span}_\C\{\hf\})^\perp
$$
be the isomorphism given by
$$
P(Z)=Z-\la Z,\hf\ra i(\hf^\perp) - \la Z,J\hf\ra Ji(\hf^\perp).
$$

\begin{proposition}\po\label{p:singsec}
The singular set of $f$ in the parametrization \eqref{e:param} is
$$
S=\{\xi\in\cl: \hat A_{i(\hfp)+\xi} \ {\rm is\ singular}\},
$$
where $\hat A=A^\hf$. The shape operator of $f$ in the direction $\hf$
restricted to $\Delta^\perp$ is
\begin{equation}\label{e:sff}
A_\hf = P(\hat A_{i(\hfp)+\xi})^{-1} P^{-1}.
\end{equation}
In particular, $S$ is also the singular set of the
submanifold itself.
\end{proposition}

\proof Take $x\in\hm$ and $\xi \in \cl(x)$. From \eqref{e:param} we see
that $df_\xi$ is the identity on $\cl(x)$. Notice that any vector
transversal to $\cl(x)$ at $\xi$ can be written as
$\psi_{*x}Z=d\psi_x(Z)$ for some $Z\in T_x\hat M$ and
$\psi\in \Gamma(\cl)$ such that $\psi(x)=\xi$. Since $\hf$ is always
normal to $f$, we have
\begin{equation}\label{e:dif}
(f_{*\xi}(\psi_{*x}Z))_{\Delta^\perp(x)}
=((f\circ\psi)_{*\xi_x}Z)_{\Delta^\perp(x)}
=((i(\hfp)+\psi)_{*x}Z)_{\Delta^\perp(x)}
=-P(\hat A_{i(\hfp)+\xi}Z),
\end{equation}
where a subspace as a subindex means to take its orthogonal projection.
For the last equality, first observe that
$$
\left((i(\hfp)+\psi)_{*x}Z\right)_{\Delta^\perp(x)}
=\left(-\hat A_{i(\hfp)+\xi}Z +\lambda_1 i(\hfp)
+ \lambda_2 Ji(\hfp) \right)_{\Delta^\perp(x)}
$$
and then use that $\hf$ is always normal to $f$ to compute the functions
$\lambda_j$, $j=1,2$. The first claim now follows from the fact that the
right-hand side of \eqref{e:dif} depends only on $Z$ and not on $\psi$.

The second part now follows since
$$
f_{*\xi}(A_\hf(\psi_{*x}Z)) = -P((\hf\circ\pi)_*(\psi_{*x}Z)) = -PZ =
P(\hat A_{i(\hfp)+\xi})^{-1}P^{-1}((f_{*\xi}(\psi_{*x}Z))_{\Delta^\perp(x)}),
$$
as we wanted.\qed
\vspace{1.5ex}

Recall that {\it first normal space\/}  of $\hf$ at $x\in \hm$
is the subspace of $N^1_\hf(x)\subseteq T_x^\perp \hat M$ spanned by the
image of the second fundamental form of $\hf$ at $x$. Equivalently,
$$
N^1_\hf(x)=\{\delta\in T_x^\perp \hat M: \hat A_{\delta}=0\}^\perp,
$$
where the orthogonal complement is taken in the normal bundle.

\begin{corollary}\po\label{c:comp}
Let $\imc f M n 1$ and $\hf\colon\hm \to \C^{n+1}$ be as in
\pref{p:antiholo}. If $M^n$ is complete, then $\hat A_{i(\hfp)}$ is
non-singular and $\cl\subseteq (N^1_\hf)^\perp$.
\end{corollary}
\proof
Assume that the conclusion does not hold. Hence, the polynomial
$$
q(z)=\det\left(\hat A_{i(\hfp)}+z\hat A_{\xi}\right)
$$
has a complex root
$u+iv$, associated to an eigenvector $U+iV\neq 0$ of the corresponding
complexified endomorphism, that is,
$$
(\hat A_{i(\hfp)}+(u+iv)\hat A_\xi)(U+iV)=0.
$$
But this is equivalent to
$$
\hat A_{i(\hfp)}U + u\hat A_\xi U - v\hat A_\xi V=0
\;\;\;\mbox{and}\;\;\;
\hat A_{i(\hfp)}V + v\hat A_\xi U + u\hat A_\xi V=0.
$$
In turn, using $J\hat A_\delta=\hat A_{J\delta}
=-\hat A_\delta J,$ we easily see that this is equivalent to
$$
\hat A_{{i(\hfp)}+(uI+vJ)\xi}(U-JV)=0\;\;\;\mbox{and}\;\;\;
\hat A_{i(\hfp)+(uI-vJ)\xi}(U+JV)=0.
$$
Since $\cl$ is holomorphic and the leaves of relative nullity are
complete, we get a contradiction with \pref{p:singsec}, because either
$U-JV$ or $U+JV$ is non-zero.
\qed
\vspace{1.5ex}
\begin{remark}\po\label{r:1leaf}
{\rm Observe that the previous result holds along each complete relative
nullity leaf of $f$, even if the submanifold is not itself complete.}
\end{remark}

As an application of \tref{t:param} we give a simple and direct proof of
Abe's cylinder theorem (\cite{abe}).

\begin{corollary}\label{c:abe}\po
Let $\imc f M n 1$ be a \hisim of a complete \K Riemannian manifold.
If the index of relative nullity satisfies $\nu \geq n-1$ everywhere, then
$f$ is an $(n-1)$-cylinder, that is, $M^n=M_1^1\times\C^{n-1}$, and
there is $f_1\colon M^1_1\to \C^2$ such that $f=f_1\times Id$ splits.
\end{corollary}
\proof
If $f$ is not totally geodesic, for which the result trivially holds, by
the hypothesis $\hf$ is an anti-holomorphic curve and, by \cref{c:comp},
we have
$$
\cl = (N^1_\hf)^\perp\;\;\;\mbox{and}\;\;\;  \spa_\C\{i(\hfp)\}=N^1_\hf.
$$
But since $\cl$ is orthogonal to the position vector $\hf$, we conclude
that the first normal space is parallel since
$0=\la \psi_*Z,\hf\ra = \la \psi_*Z,i(\hfp)\ra$ for any
$\psi\in \Gamma(\cl)$. This parallelism implies that $\hf$ reduces
codimension, that is, it is an anti-holomorphic plane curve inside some
$\C^2\subset\C^{n+1}$, and $\cl$ is the orthogonal complement of this
plane.
\qed

\section{The parametrization in $\CP^{\!N}$ } 

We show next that our parametrization in $\C^{n+1}$ can be used to obtain
a similar parametrization for holomorphic hypersurfaces of $\CP^{n+1}$.
The latter is cleaner than the former since it does not have the
restriction about the position vectors to be nowhere tangent,
and the bundle used to parametrize is the (projectivized) normal bundle
itself and not a subbundle of it.
\vspace{1.5ex}

The condition in \tref{t:param} that the position vector of a
complex submanifold $f\colon M^n\to\C^{n+p}$ is never tangent is
equivalent to the {\it cone $f_c$ over $f$} to be an immersion, where
the map $f_c\colon \C_*\times M^n\to\C^{n+p}$ is given by
$f_c(z,x)=zf(x)$. Here we denote as usual
$\C^N_*=\C^N\!\setminus\!\{0\}$. Moreover, $f$ has index of relative
nullity $\nu$ if and only if $f_c$ has index of relative nullity $\nu+1$,
and the position vector of the cone, that now is everywhere tangent,
belongs to the relative nullity.
Equivalently, the position vector of $f$ is never tangent if and only if
$f^\sharp=\hat\pi\circ f\colon M^n\to\CP^{n+p-1}$ is an immersion, where
$\hat\pi\colon \C_*^N\to\CP^{N-1}$ denotes the projection to the quotient,
and $f$ and $f^\sharp$ have the same index of relative nullity.
We conclude that to understand the submanifolds with constant relative
nullity $\nu_0$ in $\CP^N$ is equivalent to understand the cones in
$\C^{N+1}$ with $\nu\equiv\nu_0+1$.

We claim that the latter are described
as in \eqref{e:param}, but without the term $i(\hf^\perp)$.
Let $f_1\colon\,M_1^{n-1}\to\,\C^n\subset\C^{n+1}$ be the isometric
immersion obtained as the intersection of a cone
$f\colon\,M^n\to\C^{n+1}$ with constant index of relative nullity $\nu+1$
with a hyperplane, say, $\C^n=\{z_{n+1}=1\}$, so that $f_1$ has constant
index of relative nullity $\nu$ and is never tangent.
By \tref{t:param}, we have a parametrization of $f_1$ in
$\C^n$ as
$$
f_1(\xi) = i(\hf_1^\perp) + \xi, \ \ \ \xi \in \cl_1.
$$
Thus, $f_1=(\eta+\xi,1)$ in $\C^{n+1}$ where $\eta=i(\hf_1^\perp)$.
Hence, we may parametrize $f$ as
$$
f(w,\xi)=w(\eta,1)+\xi, \ \ \ \xi \in \cl_1,\ w\in\C_*.
$$
Setting $\hf=(\hf_1,-1)$, we thus parametrize $f$ as
$$
f(\xi) = \xi, \ \ \ \xi \in \cl,
$$
where $\hat f\colon \hat M^{n-\nu}\to\C^{n+2}$ is never tangent, and
$$
\cl=\cl_1\oplus{\rm span}_\C\{(\eta,1)\}=
\spa_\C\{\hat f^\perp\}^\perp\subset T^\perp_{\hat f}\hat M^{n-\nu}.
$$
This proves our claim.

 From the above description of the cones, we get for the immersion
$f^\sharp$ a parametrization over the projectivized bundle
$\P(\cl)$ of $\cl$, namely, $f^\sharp\colon \P(\cl)\to\CP^{n+1}$,
\begin{equation}\label{e:cone}
f^\sharp(\xi) = \xi, \ \ \ \xi \in \P(\cl).
\end{equation}
Now, observe that $\cl$ coincides with the normal space of $\hat f_c$,
once we identify the fibers of the normal space of $\hat f_c$ when
translated along the lines inside the cone that pass through~$0$ (we are
allowed to do this because these are lines of relative nullity of
$\hat f_c$). In other words, we have a natural identification between the
normal space of $\hat f^\sharp$ and $\cl$, and hence we can treat both as
the same fiber bundle. In particular, the corresponding complex
projectivized bundles are also identified:
$\P(T^\perp_{\hat f^\sharp}\hat M)=\P(\cl)$. These are holomorphic fiber
bundles of dimension $n$ with $\CP^\nu$ fibers. We conclude from
\eqref{e:cone} the following.

\begin{theorem}\po\label{t:paramc}
Let $\hf\colon \hm \to \CP^{n+1}$ be an anti-\hisim of a \K manifold with
vanishing relative nullity. Then, the map
$f\colon  \P(T^\perp_{\hat f}\hat M) \to \CP^{n+1}$ defined as
$$
f(\xi) = \xi,
$$
parametrizes, at regular points, a \h \ \K hypersurface with constant
index of relative nullity $\nu$.
Conversely, any such hypersurface can be parametrized this way.
\end{theorem}

We point out that the holomorphicity hypothesis in the converse is
redundant when the submanifold has relative nullity. It was shown in
\cite{dr} that any isometric immersion of a \K manifold into $\CP^N$
with positive index of relative nullity must be holomorphic.

\begin{remark}\po\label{r:fin}
{\rm
Taking the Gauss map is an involution
on $\bar{\cal H}_+\cup\bar{\cal H}_-$ that swaps $\bar{\cal H}_+$ with
$\bar{\cal H}_-$, where $\bar{\cal H}_+$ and $\bar{\cal H}_{-}$ denote
the sets of holomorphic and anti-holomorphic hypersurfaces of $\CP^{n+1}$
with vanishing relative nullity. As opposed to the $\C^N$ case, here
there is no restriction on the position vectors; see \rref{r:invol}.
}
\end{remark}

{\renewcommand{\baselinestretch}{1}
\hspace*{-20ex}\begin{tabbing}
\indent \= IMPA -- Estrada Dona Castorina, 110 \\
\> 22460-320 --- Rio de Janeiro --- Brazil \\
\> \ marcos@impa.br \ \ --- \ \  luis@impa.br
\end{tabbing}
}
\end{document}